\documentclass[12pt]{amsart}

\usepackage[utf8]{inputenc}
\usepackage[T1]{fontenc}

\usepackage{amssymb}
\usepackage{amsmath}
\usepackage{amsfonts}
\usepackage{latexsym}
\usepackage{amsthm}
\usepackage{mathrsfs} 
\usepackage[colorlinks, urlcolor=blue, citecolor=black, linkcolor=black ,breaklinks]{hyperref} 

\usepackage{enumerate}

\newtheorem{Definition}{Definition}[section]
\newtheorem{Theorem}[Definition]{Theorem}
\newtheorem{Lemma}[Definition]{Lemma}
\newtheorem{Corollary}[Definition]{Corollary}
\newtheorem{Proposition}[Definition]{Proposition}

\newcommand{\Remark}[1][ ]{{\bf Remark #1: }}

\def\R{\mathbb{R}}

\def\N{\mathbb{N}}
\def\F{\mathbb{F}}

\newcommand{\Ld}[2][D]{\ensuremath{\mathcal{L}(#2\mathcal{#1})}} 
\newcommand{\D}[1][D]{\ensuremath{\mathcal{#1}}}

\begin{document}

\title[On the tensor rank of multiplication in $\F_{2^n}$]{On the tensor rank of multiplication in any extension of $\F_2$}
\author{St\'ephane Ballet}
\address{Institut de Math\'{e}matiques de Luminy\\ case 930, F13288 Marseille cedex 9\\ France}
\address{eRISCS Groupe de Recherche "Informatique des Syst\`{e}mes Communicants S\'{e}curis\'{e}s"}
\email{stephane.ballet@univmed.fr}
\author{Julia Pieltant}
\address{Institut de Math\'{e}matiques de Luminy\\ case 930, F13288 Marseille cedex 9\\ France}
\address{eRISCS Groupe de Recherche "Informatique des Syst\`{e}mes Communicants S\'{e}curis\'{e}s"}
\email{julia.pieltant@univmed.fr}
\date{\today}
\keywords{Algebraic function fields, tower of function fields, tensor rank, algorithm, finite fields}
\subjclass[2000]{Primary 14H05; Secondaries 11Y16, 12E20}

\begin{abstract}
In this paper, we obtain new bounds for the tensor rank of multiplication in any extension of $\F_2$. In particular, it also enables us to obtain the best known asymptotic bound. In this aim, we use the generalized algorithm of type Chudnovsky with derivative evaluations on places of degree one, two and four applied on  the descent over $\F_2$ of a Garcia-Stichtenoth tower of algebraic function fields defined over $\F_{2^4}$.

\end{abstract}

\maketitle


\section{Introduction}

\subsection{Tensor rank of multiplication}

Let $\F_q$ be a finite field with $q$ elements where $q$ is a prime power and let $\F_{q^n}$ be a $\F_q$-extension of degree $n$. We denote by $m$ the ordinary 
multiplication  in the $\F_q$-vector space $\F_{q^n}$. The multiplication $m$ is a bilinear map from ${\F_{q^n} \times \F_{q^n}}$ into $\F_{q^n}$, thus it corresponds to a linear map $M$ from the tensor product ${\F_{q^n} \bigotimes \F_{q^n}}$ over $\F_q$ into $\F_{q^n}$. One can also represent $M$ by a tensor ${t_M \in \F_{q^n}^*\bigotimes \F_{q^n}^* \bigotimes \F_{q^n}}$ where $\F_{q^n}^*$ denotes the dual of $\F_{q^n}$ over $\F_q$. Hence the product of two elements $x$ and $y$ of $\F_{q^n}$ is the convolution of this tensor with ${x \otimes y \in \F_{q^n} \bigotimes \F_{q^n}}$. If 
\begin{eqnarray}
t_M = \sum_{l=1}^{\lambda} a_l \otimes b_l \otimes c_l
\end{eqnarray}
where ${a_l \in  \F_{q^n}^*}$, ${b_l \in  \F_{q^n}^*}$, ${c_l \in  \F_{q^n}}$, then 
\begin{eqnarray}
x\cdot y=\sum_{l=1}^{\lambda} a_l(x)b_l(y)c_l.
\end{eqnarray}

Every expression (2) is called a bilinear multiplication algorithm $\mathcal U$. The integer $\lambda$ is called the multiplicative complexity  ${\mu({\mathcal U})}$ of $\mathcal U$.\\
Let us set
$$
\mu_{q}(n)= \min_\mathcal{U} \mu(\mathcal{U}),
$$
where $\mathcal U$ is running over all bilinear multiplication algorithms in $\F_{q^n}$ over $\F_q$.\\
Then $\mu_{q}(n)$ corresponds to the minimum possible number of summands in any tensor decomposition of type (1), which is the rank of 
the tensor of multiplication in $\F_{q^n}$ over $\F_q$. The tensor rank $\mu_{q}(n)$ is also called the bilinear complexity of multiplication in $\F_{q^n}$ over $\F_q$.

\subsection{Notations} 

Let $F/\F_q$ be an algebraic function field of one variable of genus $g$, with constant field $\F_q$, associated to a curve $X$ defined over $\F_q$. 
For any place $P$ we define $F_P$ to be the residue class \linebreak[4]field of $P$ and $\mathcal{O}_P$ its valuation ring. Every element ${t \in P}$ such that ${P = t \mathcal{O}_P}$ is called a local parameter for $P$ and we denote by ${v_P}$ a discrete valuation associated to the place $P$ in $F/\F_q$. Recall that this valuation does not depend on the choice of the local parameter. \linebreak[4] If $\D$ is a divisor then 
${\Ld{}=\{f \in F/\F_q ; \D + (f) \geq 0\}\cup \{0\}}$ is a vector space over $\F_q$ whose dimension $\dim {\D}$ is given by the 
Riemann-Roch Theorem. 
The degree of a divisor ${\D}=\sum_P a_P P$ is defined by ${\deg \D=\sum_P a_P \deg P}$ where $\deg P$ is the dimension of $F_P$ over $\F_q$. 
The order of a divisor ${\D=\sum_P a_P P}$ in $P$ is the integer $a_P$ denoted by ${\mathrm{ord}_P \, \D}$. The support of a divisor $\D$ is the set ${\mathrm{supp} \, \D}$ of the places $P$ such that ${\mathrm{ord}_P \, \D \neq 0}$. Let ${x \in F\backslash \{0\}}$, we denote by ${(x) := \sum_P v_P(x) P}$ where $P$ is running over all places in $F/\F_q$, the principal divisor of $x$. Two divisors $\D$ and $\D'$ are said to be equivalents if ${\D=\D'+(x)}$ for an element ${x \in F\backslash \{0\}}$.

\subsection{Known results} 

\subsubsection{General results}

The bilinear complexity $\mu_q(n)$ of the multiplication in the $n$-degree extension of a finite field $\F_q$  with $q$ elements is known for certain values of $n$.  In particular, S. Winograd \cite{wino3} and H. De Groote \cite{groo} have shown that this complexity is ${\geq 2n-1}$, with equality holding if and only if ${n \leq \frac{1}{2}q+1}$. 
Using the principle of the D.V. and G.V. Chudnovsky algorithm \cite{chch} applied to elliptic curves, M.A. Shokrollahi has shown in \cite{shok} that the bilinear complexity of multiplication is equal to $2n$ for ${\frac{1}{2}q +1< n < \frac{1}{2}(q+1+{\epsilon (q) })}$ where $\epsilon$ is the function defined by:
$$
\epsilon (q) = \left \{
	\begin{array}{l}
 		 \mbox{greatest integer} \leq 2{\sqrt q} \mbox{ prime to $q$, if $q$ is not a perfect square} \\
  		2{\sqrt q}\mbox{, if $q$ is a perfect square.}
	\end{array} \right .
$$
Moreover, U. Baum and M.A. Shokrollahi have succeeded in \cite{bash} to construct effective optimal algorithms of type Chudnovsky in the elliptic case. 

Recently in \cite{ball1}, \cite{ball3}, \cite{baro1}, \cite{balbro}, \cite{balb}, \cite{bach} and \cite{ball5} the study made by M.A. Shokrollahi has been  generalized to algebraic function fields of genus~$g$. 

In particular, from the existence of towers of algebraic function fields satisfying good properties, it was proved:

\begin{Theorem} \label{theobornes}
Let $q=p^r$ be a power of the prime $p$. The tensor rank $\mu_q(n)$ of multiplication in any finite field $\F_{q^n}$ is linear with respect to the extension degree, more precisely:
$$
\mu_q(n) \leq C_q n
$$
where $C_q$ is the constant defined by: 
$$
C_q = \left \{
	\begin{array}{lll}
		\hbox{if } q=2 & \hbox{then } 54 &   \hbox{\cite{ball1}} \cr
		\hbox{if } q=3 &  \mbox{then } 27 & \hbox{\cite{ball1}} \cr
		\hbox{else if } q=p \geq 5 &  \mbox{then }   3(1+ \frac{4}{q-3}) & \hbox{\cite{bach}} \\
		\hbox{else if } q=p^2 \geq 25 & \hbox{then }   2(1+\frac{2}{\sqrt{q}-3}) & \hbox{\cite{bach}} \\
		\hbox{else if } q=p^{2k} \geq 16 & \hbox{then } 2(1+\frac{p}{\sqrt{q}-3}) & \hbox{\cite{ball3}} \\
		\hbox{else if } q \geq 16 & \hbox{then }  3(1+\frac{2p}{q-3}) & \hbox{\cite{baro1}, \cite{balbro} and \cite{balb}}\\
		\hbox{else if } q > 3 & \hbox{then }  6(1+\frac{p}{q-3}) & \hbox{\cite{ball3}}.
	\end{array}\right .
$$
\end{Theorem}

In order to obtain these good estimates for the constant $C_q$, S. Ballet has given in \cite{ball1} some easy to verify conditions allowing the use of the D.V. and G.V. Chudnovsky  algorithm. Then S. Ballet and R. Rolland have generalized in \cite{baro1} the algorithm using places of degree $1$ and~$2$.

Let us present the last version of this algorithm, which is a generalization of the algorithm of type Chudnovsky introduced by N. Arnaud in \cite{arna1} and M. Cenk and F. \"Ozbudak in \cite{ceoz}. This generalization uses several coefficients in the local expansion at each place $P_i$ instead of just the first one. Due to the way to obtain the local expansion of a product from the local expansion of each term, the bound for the bilinear complexity involves the complexity notion $\widehat{M_q}(u)$ introduced by M. Cenk and F. \"Ozbudak in \cite{ceoz} and defined as follows:
\begin{Definition}
We denote by $\widehat{M_q}(u)$ the minimum number of multiplications needed in $\F_q$ in order to obtain coefficients of the product of two arbitrary $u$-term polynomials modulo $x^u$ in $\F_q[x]$.
\end{Definition}
For instance, we know that for all prime powers $q$, we have $\widehat{M_q}(2) \leq 3$ by \cite{ceoz2}.\\

Now we introduce the generalized algorithm of type Chudnovsky described in~\cite{ceoz}.

\begin{Theorem} \label{theo_evalder}
Let \\
\vspace{.1em}
$\bullet$ $q$ be a prime power,\\
\vspace{.1em}
$\bullet$ $F/\F_q$ be an algebraic function field,\\
\vspace{.1em}
$\bullet$ $Q$ be a degree $n$ place of $F/\F_q$,\\
\vspace{.1em}
$\bullet$ $\D$ be a divisor of $F/\F_q$,\\
\vspace{.1em}
$\bullet$ ${\mathcal P}=\{P_1,\ldots,P_N\}$ be a set of $N$ places of arbitrary degree,\\
\vspace{.1em}
$\bullet$ $u_1,\ldots,u_N$ be positive integers.\\
We suppose that $Q$ and all the places in $\mathcal P$ are not in the support of $\D$ and that:
\begin{enumerate}[a)]
	\item the application
	$$
	Ev_Q:  \left \{
	\begin{array}{ccl}
	\Ld{} & \rightarrow & \F_{q^n}\simeq F_Q\\
	f & \longmapsto & f(Q)
	\end{array} \right.
	$$ 
	is onto,
	\item the application
	$$
	Ev_{\mathcal P} :  \left \{
	\begin{array}{ccl}
	\Ld{2} & \longrightarrow & \left(\F_{q^{\deg P_1}}\right)^{u_1} \times \left(\F_{q^{\deg P_2}}\right)^{u_2} \times \cdots \times \left(\F_{q^{\deg P_N}}\right)^{u_N} \\
	f & \longmapsto & \big(\varphi_1(f), \varphi_2(f), \ldots, \varphi_N(f)\big)
	\end{array} \right.
	$$
	is injective, where the application $\varphi_i$ is defined by
	$$
	\varphi_i : \left \{
	\begin{array}{ccl}
	\Ld{2} & \longrightarrow & \left(\F_{q^{\deg P_i}}\right)^{u_i} \\
          f & \longmapsto & \left(f(P_i), f'(P_i), \ldots, f^{(u_i-1)}(P_i)\right)
	\end{array} \right.
 	$$
	with $f = f(P_i) + f'(P_i)t_i + f''(P_i)t_i^2+ \ldots + f^{(k)}(P_i)t_i^k + \ldots $, the local expansion at $P_i$ of $f$ in ${\Ld{2}}$, with respect to the local parameter~$t_i$. Note that we set ${f^{(0)} =f}$.

\end{enumerate}
Then 
$$
\mu_q(n) \leq \displaystyle \sum_{i=1}^N \mu_q(\deg P_i) \widehat{M}_{q^{\deg P_i}}(u_i).
$$
\end{Theorem}

First of all, note that we can define the application $Ev_Q$ since $Q$ is not in the support of $\D$. Indeed, for such a place $Q$, we have ${\Ld{} \subseteq \mathcal{O}_Q}$, so $Ev_Q$ is the restriction of the projection ${\pi : \mathcal{O}_Q \rightarrow F_Q}$. Moreover, the application $Ev_\mathcal{P}$ can be define since ${\Ld{2} \subseteq \mathcal{O}_{P_i}}$ for all intergers ${i \in \{1, \ldots, n\}}$, so the local expansion of ${f\in \Ld{2}}$ at any place $P_i \in \mathcal P$ exists from \cite{nixi2} (1.4). Indeed, this follows from the fact that the intersection ${\mathcal{P} \cap \mathrm{supp} \, \D}$ is empty, so ${v_{P_i}(f) \geq 0}$ and the coefficients of the local expansion of $f$ at $P_i$ can be define inductively.\\

Let us remark that the algorithm given in \cite{chch} by D.V. and G.V. Chudnovsky is the case $\deg P_i=1$ and $u_i=1$ for $i=1, \ldots, N$. The generalization introduced here is usefull: it allows us to use certain places many times, thus less places are necessary to get the injectivity of $Ev_\mathcal{P}$.

In particular, we have the following results, obtained by N. Arnaud in \cite{arna1}. 

\begin{Corollary} \label{theo_deg12evalder}
Let \\
\vspace{.1em}
$\bullet$ $q$ be a prime power,\\
\vspace{.1em}
$\bullet$ $F/\F_q$ be an algebraic function field,\\
\vspace{.1em}
$\bullet$ $Q$ be a degree $n$ place of $F/\F_q$,\\ 
\vspace{.1em}
$\bullet$ $\D$ be a divisor of $F/\F_q$, \\
\vspace{.1em}
$\bullet$ ${\mathcal P}=\{P_1,\ldots,P_{N_1},P_{N_1+1},\ldots,P_{N_1+N_2}\}$ be a set of $N_1$ places of degree\\ \indent one and $N_2$ places of degree two,\\
\vspace{.1em}
$\bullet$ ${0 \leq l_1 \leq N_1}$ and ${0 \leq l_2 \leq N_2}$ be two integers.\\
We suppose that $Q$ and all the places in $\mathcal P$ are not in the support of $\D$ and that:
\begin{enumerate}[a)]
	\item the application
	$$
	Ev_Q: \Ld{} \rightarrow \F_{q^n}\simeq F_Q$$
	is onto,
	\item the application
	 $$
	 Ev_{\mathcal P}: \left \{
	\begin{array}{ccl}
   	\Ld{2} & \rightarrow & \F_{q}^{N_1} \times \F_{q}^{l_1}\times \F_{q^2}^{N_2} \times \F_{q^2}^{l_2}  \\
    	f & \mapsto & \big(f(P_1),\ldots,f(P_{N_1}),f'(P_1),\ldots,f'(P_{l_1}),\\
	  &  & \ f(P_{N_1+1}),\ldots,f(P_{N_1+N_2}),f'(P_{N_1+1}),\ldots,f'(P_{N_1+l_2})\big)
	\end{array} \right .
	$$
is injective.
\end{enumerate}
Then  
$$
\mu_q(n)\leq  N_1 + 2l_1 + 3N_2 + 6l_2.
$$
\end{Corollary}

\begin{Proof}
Up to reindex the places, the result follows from Theorem~\ref{theo_evalder} applied with $N=N_1+N_2$,  ${\deg P_i=1}$ for $i=1, \ldots, N_1$, ${\deg P_i = 2}$ for $i=N_1+1, \ldots, N$, and 
$$u_i = \left \{
	\begin{array}{ll}
	2\mbox{,} & \mbox{if }1\leq i \leq l_1 \mbox{ or } N_1+1 \leq i \leq N_1+l_2,\\
	1\mbox{,} & \mbox{else.}
	\end{array} \right. $$

Recall that for all prime powers $q$, we have $\mu_q(2)=3$ and $\widehat{M_q}(2) \leq 3$.\\
Then applying Theorem \ref{theo_evalder}, we get:
\begin{eqnarray*}
	\mu_q(n) & \leq & \sum_{i=1}^{l_1}\mu_q(1)\widehat{M_q}(2) + \sum_{i=l_1+1}^{N_1}\mu_q(1)\widehat{M_q}(1) + \sum_{i=N_1+1}^{N_1+l_2}\mu_q(2)\widehat{M_{q^2}}(2)\\
			& & \  + \sum_{i=N_1+l_2+1}^{N}\mu_q(2)\widehat{M_{q^2}}(1)\\
			& \leq & 3l_1 + N_1 - l_1 +9l_2 + 3(N_2-l_2)\\
			& = & N_1 +2l_1+3N_2 +6l_2.
\end{eqnarray*}\qed\\
\end{Proof}

In particular, from the last corollary applied on Garcia-Stichtenoth towers, N. Arnaud obtained the two following bounds.

\begin{Theorem}\label{theo_arnaud1}
Let ${q=p^r \geq 4}$ be a prime power. Then
	\begin{equation*}
	\mbox{(i) \ \ \ \ \ }\mu_{q^2}(n) \leq 2 \left(1 + \frac{p}{q-3 + (p-1)\left(1- \frac{1}{q+1}\right)} \right)n,
	\end{equation*}
	\begin{equation*}
	\mbox{(ii) \ \ \ \ \ }\mu_{q}(n) \leq 3 \left(1 + \frac{2p}{q-3 + 2(p-1)\left(1- \frac{1}{q+1}\right)} \right)n.
	\end{equation*}

\end{Theorem}

\subsubsection{Asymptotic bounds for the extensions of $\F_2$}

From the asymptotical point of view, let us recall that I. Shparlinski, M. Tsfasman and S. Vladut have given in \cite{shtsvl} many interesting remarks on the algorithm of D.V. and G.V. Chudnovsky. In particular, they considered the following asymptotic bounds for the bilinear complexity

$$
M_q = \limsup_{k \rightarrow \infty} \frac{\mu_q(k)}{k}
$$
and
$$
m_q = \liminf_{k \rightarrow \infty} \frac{\mu_q(k)}{k}.
$$
In \cite{shtsvl}, they claim that $M_2 \leq 27$ but I. Cascudo, R. Cramer and C. Xing recently noticed in \cite{xing} a gap in the proof about $M_2 \leq 27$ established by I. Shparlinsky, M. Tsfasman and S. Vladut: \textit{the mistake in \cite{shtsvl}  from 1992 is in the proof of their Lemma 3.3, page 161, the paragraph following formulas about the degrees of the divisor. It reads: "}\textsl{Thus the  number of linear equivalence classes of degree a for which either Condition $\alpha$ or Condition $\beta$ fails is at most $D_{b'} + D_b$.}\textit{" This is incorrect; $D_b$ should be multiplied by the torsion. Hence the proof of their asympotic bound is incorrect.} \\
Consequently, we can consider this bound as not proved. Anyway, it is possible to obtain easily a better bound for $M_2$ from one of the bounds of N. Arnaud. Indeed, by using Bound (ii) of Theorem \ref{theo_arnaud1}.

we obtain:

\begin{Proposition}\label{theo_M2_2285}
$$
M_2 \leq \frac{297}{13} \approx 22.85.
$$
\end{Proposition}

\begin{Proof} 
For all $m \geq 1$, we have
$$
\mu_q(n) \leq \mu_q(mn) \leq \mu_q(m) \cdot \mu_{q^m}(n).
$$
Thus for ${q=2}$ and ${m=2}$ we get ${\mu_2(n) \leq \mu_2(2) \cdot \mu_{4}(n)}$. Remembering that $\mu_2(2)=3$ and applying Bound (ii) of Theorem \ref{theo_arnaud1}, we obtain 
$$
\mu_2(n) \leq 3\cdot 3 \left(1+ \frac{4}{1 + 2 \left(1-\frac{1}{5}\right)}\right)n = \frac{297}{13}n.
$$
\qed
\end{Proof}

\Remark Using Bound (i) from Theorem \ref{theo_arnaud1}, we obtain $M_2 \leq 38$. Indeed, for all $m \geq 1$, we have
$$
\mu_q(n) \leq \mu_q(mn) \leq \mu_q(m) \cdot \mu_{q^m}(n).
$$
Thus for ${q=2}$ and ${m=4}$ we get ${\mu_2(n) \leq \mu_2(4) \cdot \mu_{16}(n)}$. Remembering that $\mu_2(4)\leq9$ and applying Bound (i) of Theorem \ref{theo_arnaud1}, we obtain 
$$
\mu_2(n) \leq 9\cdot 2 \left(1+ \frac{2}{2-\frac{1}{5}}\right)n = 38n.
$$

\subsection{New results established in this paper}  Our main result concerns an improvement of the asymptotic bound for the tensor rank of multiplication in any extension of $\F_2$. More precisely, we prove that:
$$
M_2 \leq \frac{477}{26} \approx 18.35.
$$
This result comes from a new bound for the tensor rank of multiplication in any extension of $\F_2$ that we also obtain in this paper, namely:
$$
\mu_2(n) \leq \frac{477}{26}n +\frac{45}{2}.
$$
In Section 2, we recall some results about a modified Garcia-Stichtenoth's tower  \cite{gast} studied in \cite{ball3}, \cite{baro1}, \cite{balbro} and \cite{bach}. Specially, we present the descent of  the definition field of this Garcia-Stichtenoth's tower on the field $\F_2$  obtained in \cite{baro4} and study some of its properties which will be useful in Section 4. In Section 3, we specialize the generalized algorithm of type Chudnovsky by using places of degree one, two and four with derivative evaluations. In order to obtain new bounds for the bilinear complexity, we apply this specialized algorithm to suitable steps of the tower presented in Section 2. In particular, in Section 4 these new bounds lead to an improvement of known results on the asymptotic tensor rank of multiplication in the extensions of $\F_2$.

\section{A good  sequence of function fields defined over $\F_2$} \label{asympexactsequence}

In this section, we present a sequence of algebraic function fields defined over $\F_2$ constructed and studied in \cite{baro4}, which will be used to obtain the new bounds of the tensor rank of multiplication in the extensions of $\F_2$.

\subsection{Definition of  Garcia-Stichtenoth's towers} \label{sectiondefinition}
 
First, we  present a modified Garcia-Stichtenoth's tower (cf. \cite{gast}, \cite{ball3}, \cite{baro1}) having good properties. Let  us consider a finite field $\F_{q^2}$ with ${q=p^r}$, for $p$ a prime number and $r$ an integer. Let us consider the Garcia-Stichtenoth's elementary abelian tower $T_0$ over $\F_{q^2}$ constructed in \cite{gast} and defined by the sequence ${(F_1, F_2,\ldots)}$ where 
$$
F_{k+1}:=F_{k}(z_{k+1})
$$ 
and $z_{k+1}$ satisfies the equation: 
$$
z_{k+1}^q+z_{k+1}=x_k^{q+1}
$$ 
with 
$$
x_k:=z_k/x_{k-1} \mbox{ in } F_k \mbox{ (for } k \geq1 \mbox{).}
$$
Moreover ${F_1:=\F_{q^2}(x_0)}$ is the rational function field over $\F_{q^2}$ and  $F_2$ the Hermitian function field over $\F_{q^2}$. Let us denote by $g_k$ the genus of $F_k$ in $T_0/\F_{q^2}$, we recall the following \textsl{formulae}:
\begin{equation}\label{genregs}
g_k = \left \{ \begin{array}{ll}
		q^k+q^{k-1}-q^\frac{k+1}{2} - 2q^\frac{k-1}{2}+1 & \mbox{if } k \equiv 1 \mod 2,\\
		q^k+q^{k-1}-\frac{1}{2}q^{\frac{k}{2}+1} - \frac{3}{2}q^{\frac{k}{2}}-q^{\frac{k}{2}-1} +1& \mbox{if } k \equiv 0 \mod 2.
		\end{array} \right .
\end{equation}
If $r>1$, we consider the completed Garcia-Stichtenoth's tower 
$$
T_1/\F_{q^2}= F_{1,0}\subseteq F_{1,1}\subseteq \ldots\subseteq F_{1,r} = F_{2,0}\subseteq F_{2,1} \subseteq\ldots\subseteq F_{2,r} = F_{3,0} \subseteq \ldots
$$ 
considered in \cite{ball3} such that ${F_k \subseteq F_{k,s} \subseteq F_{k+1}}$ for any integer $s$ such that ${s=0,\ldots,r}$, with ${F_{k,0}=F_k}$ and ${F_{k,r}=F_{k+1}}$. 
Let us denote by $g_{k,s}$ the genus of $F_{k,s}/\F_{q^2}$ in $T_1/\F_{q^2}$ and by $N_i(F_{k,s}/\F_{q^2})$ the number of places of degree$~i$ of $F_{k,s}/\F_{q^2}$ in $T_1/\F_{q^2}$. Recall that each extension $F_{k,s}/F_k$ is Galois of degree $p^s$ wih full constant field $\F_{q^2}$. Moreover, we know by \cite{balbro} that the descent of the definition field of the tower $T_1/\F_{q^2}$ from $\F_{q^2}$ to  $\F_{q}$ is possible. More precisely,  there exists a tower $T_2/\F_q$ defined over $\F_{q}$ given by a sequence:
$$
T_2/\F_q=G_{1,0} \subseteq G_{1,1} \subseteq G_{1,2} = G_{2,0}\subseteq G_{2,1}\subseteq G_{2,2} =G_{3,0} \subseteq \ldots
$$
defined over the constant field $\F_q$ and related to the tower $T_1/\F_{q^2}$ by
$$
F_{k,s}=\F_{q^2}G_{k,s} \mbox{ for all $k$ and $s$,}
$$
namely $F_{k,s}/\F_{q^2}$ is the constant field extension of $G_{k,s}/\F_q$.

\subsection{Descent of the definition field of a Garcia-Stichtenoth's tower on the field $\F_2$} \label{sectiondescent}

Now, we are interested to search the descent of the definition field of the tower $T_1/\F_{q^2}$ from $\F_{q^2}$ to  $\F_{p}$ if it is possible. In fact, one can not establish a general result but one can prove that it is possible in the case of characteristic $2$ which is given by the following result obtained in \cite{baro4}. Note that in order to simplify the presentation, we are going to set the results by using the variable $p$ and to give the proofs to be self-contained.

\begin{Proposition}\label{proptour}
Let $p=2$. If $q=p^2$, the descent of the definition field of the tower $T_1/\F_{q^2}$ from $\F_{q^2}$ to  $\F_{p}$ is possible. More precisely,  there exists a tower $T_3/\F_p$ defined over $\F_{p}$ given by a sequence:
$$
T_3/\F_p=H_{1,0} \subseteq H_{1,1} \subseteq H_{1,2} = H_{2,0}\subseteq H_{2,1}\subseteq H_{2,2} = H_{3,0} \subseteq \ldots
$$
defined over the constant field $\F_p$ and related to the towers $T_1/\F_{q^2}$ and $T_2/\F_q$ by
$$
F_{k,s}=\F_{q^2}H_{k,s} \mbox{ for all $k$ and $s=0,1,2$},
$$
$$
G_{k,s}=\F_{q}H_{k,s} \mbox{ for all $k$ and $s=0,1,2$},
$$
namely $F_{k,s}/\F_{q^2}$ is the constant field extension of $G_{k,s}/\F_q$ and $H_{k,s}/\F_p$ and $G_{k,s}/\F_q$ is the constant field extension of $H_{k,s}/\F_p$. 
\end{Proposition}

\begin{Proof}
In the proof, we use $p=2$. 
Let $x_1$ be a transcendent element over $\F_2$ and let us set
$$H_1=\F_2(x_1), G_1=\F_4(x_1), F_1=\F_{16}(x_1).$$
We define recursively for $k \geq 1$
\begin{enumerate}
	\item $z_{k+1}$ such that  $z_{k+1}^4+z_{k+1}=x_k^5$,
 	\item $t_{k+1}$ such that $t_{k+1}^2+t_{k+1} = x_k^5$\\
(or alternatively $t_{k+1}=z_{k+1}(z_{k+1}+1)$),
	\item $x_{k+1}=z_{k+1}/x_k$,
	\item $H_{k,1}=H_{k,0}(t_{k+1})=H_k(t_{k+1})$,
$H_{k+1,0}= H_{k+1}=H_k(z_{k+1})$,
$G_{k,1}=G_{k,0}(t_{k+1})=G_k(t_{k+1})$,
$G_{k+1,0}= G_{k+1}=G_k(z_{k+1})$,
$F_{k,1}=F_{k,0}(t_{k+1})=F_k(t_{k+1})$,
$F_{k+1,0}= F_{k+1}=F_k(z_{k+1})$.
\end{enumerate}
By \cite{balbro}, the tower $T_1=(F_{k,i})_{k \geq 1, i=0,1}$ is the densified Garcia-Stichtenoth's tower
over $\F_{16}$ and the two other towers $T_2$ and $T_3$ 
are respectively the descent of $T_1$
over $\F_4$ and over $\F_2$.  
\qed\\
\end{Proof}

Now, we recall different properties concerning the tower $T_3/\F_2$.

\begin{Proposition}\label{subfield}
Let $q=p^2=4$. For any integers ${k\geq1}$ and ${s \in \{0,1,2\}}$, the algebraic function field $H_{k,s}/\F_{p}$ in the tower $T_3/\F_p$ has a genus ${g(H_{k,s}/\F_p)=g_{k,s}}$ with ${N_{1}(H_{k,s}/\F_p)}$ places of degree one, ${N_{2}(H_{k,s}/\F_p)}$ places of degree two and ${N_{4}(H_{k,s}/\F_p)}$ places of degree $4$ such that:
\begin{enumerate}[1)]
	\item ${H_k/\F_p \subseteq H_{k,s} /\F_p \subseteq H_{k+1}/\F_p}$ with ${H_{k,0}=H_k}$ and ${H_{k,2}=H_{k+1}}$,
	\item ${g(H_{k,s}/\F_p)\leq \frac{g(H_{k+1}/\F_p)}{p^{2-s}}+1}$ with ${g(H_{k+1}/\F_p)=g_{k+1}\leq q^{k+1}+q^{k}}$,   
	\item ${N_{1}(H_{k,s}/\F_p)+2N_{2}(H_{k,s}/\F_p) +4N_{4}(H_{k,s}/\F_p)\geq (q^2-1)q^{k-1}p^{s}}$,

\end{enumerate}
\end{Proposition}

\begin{Proof}
The property $1)$ follows directly from Proposition \ref{proptour}. Moreover, by Theorem 2.2 in \cite{ball3}, we have 
$g(F_{k,s})\leq {g(F_{k+1}) \over p^{2-s}}+1$ with $g(F_{k+1})=g_{k+1}\leq q^{k+1}+q^{k}$ . Then, 
as the algebraic function field $F_{k,s}$ is a constant field extension of $H_{k,s}$, for any integers $k$ and $s$ 
the algebraic function fields $F_{k,s}$ and $H_{k,s}$ have the same genus. So, the inequality satisfied by the 
genus $g(F_{k,s})$ is also true for the genus $g(H_{k,s})$. Moreover, the number of places of degree one 
$N_1(F_{k,s}/\F_{q^2})$ of $F_{k,s}/\F_{q^2}$ is such that $N_1(F_{k,s}/\F_{q^2})\geq (q^2-1)q^{k-1}p^{s}$. Then, as 
the algebraic function field $F_{k,s}$ is a constant field extension of $H_{k,s}$ of degree $4$, it is clear that for any
integers $k$ and $s$, we have $${N_{1}(H_{k,s}/\F_p)+2N_{2}(H_{k,s}/\F_p)+4N_{4}(H_{k,s}/\F_p)\geq (q^2-1)q^{k-1}p^{s}}.\qed$$

\end{Proof}

\subsection{Some preliminary results} \label{sectionusefull}
Here we gather some results about genus and number of places of each step of the tower $T_3/\F_2$ defined in Section \ref{sectiondescent}. These results will allow us to determine a suitable step of the tower to apply the algorithm on. In order to simplify the presentation, we still use the variables $p$ and $q$.

\begin{Lemma}\label{lemme_genre}
Let ${q=p^2=4}$. We have the following bounds for the genus of each step of the tower $T_3/\F_p$:
\begin{enumerate}[i)]
	\item $g_k> q^k$ for all ${k\geq 4}$,
	\item $g_k \leq q^{k-1}(q+1) - \sqrt{q}q^\frac{k}{2}$,
	\item $g_{k,s} \leq q^{k-1}(q+1)p^s$ for all ${k\geq 1}$, $s=0,1,2$,
	\item $g_{k,s} \leq \frac{q^k(q+1)-q^\frac{k}{2}(q-1)}{p^{2-s}}$ for all $k\geq 2$, $s=0,1,2$.
\end{enumerate}
\end{Lemma}

\begin{Proof}
\textit{i)} According to Formula (\ref{genregs}) recalled in Section \ref{sectiondefinition}, we know that if ${k \equiv 1 \mod 2}$, then 
$$
g_k = q^k+q^{k-1}-q^\frac{k+1}{2} - 2q^\frac{k-1}{2}+1 = q^k+q^\frac{k-1}{2}(q^\frac{k-1}{2} - q - 2) +1.
$$
Since ${q=4}$ and ${k \geq 4}$, we have ${q^\frac{k-1}{2} - q - 2 >0}$, thus ${g_k>q^k}$.\\
Else if ${k \equiv 0 \mod 2}$, then 
$$
g_k = q^k+q^{k-1}-\frac{1}{2}q^{\frac{k}{2}+1} - \frac{3}{2}q^{\frac{k}{2}}-q^{\frac{k}{2}-1} +1 = q^k+q^{\frac{k}{2}-1}(q^\frac{k}{2}-\frac{1}{2}q^{2} - \frac{3}{2}q-1)+1.
$$
Since ${q=4}$ and ${k\geq 4}$, we have ${q^\frac{k}{2}-\frac{1}{2}q^{2} - \frac{3}{2}q-1>0}$, thus ${g_k>q^k}$.

\textit{ii)} It follows from Formula (\ref{genregs}) since for all $k\geq 1$ we have ${2q^\frac{k-1}{2} \geq 1}$ which works out for odd $k$ cases and ${\frac{3}{2}q^\frac{k}{2}+q^{\frac{k}{2}-1}\geq 1}$ which works out for even $k$ cases. Recall that ${\frac{1}{2}q=\sqrt{q}}$ here.

\textit{iii)} If ${s=2}$, then according to Proposition \ref{subfield}, we have 
$$
g_{k,s} = g_{k+1}\leq q^{k+1}+q^{k} = q^{k-1}(q+1)p^2.
$$
Else, ${s<2}$ and Proposition \ref{subfield} says that ${g_{k,s} \leq \frac{g_{k+1}}{p^{2-s}}+1}$. Moreover, since ${q^\frac{k+2}{2}\geq q}$ and ${\frac{1}{2}q^{\frac{k+1}{2}+1}\geq q}$, we obtain ${g_{k+1}\leq q^{k+1} + q^k - q + 1}$ from Formula (\ref{genregs}). Thus, we get 
\begin{eqnarray*}
g_{k,s} & \leq & \frac{q^{k+1} + q^k - q + 1}{p^{2-s}} +1\\
	    &  = & q^{k-1}(q+1)p^s - p^s + p^{s-2} + 1\\
	    & \leq & q^{k-1}(q+1)p^s + p^{s-2}\\
	    & \leq & q^{k-1}(q+1)p^s \  \mbox{ since ${0 \leq p^{s-2} <1}$ and ${g_{k,s} \in \N}$}.
\end{eqnarray*}

\textit{iv)} It follows from ii) since Proposition \ref{subfield} gives ${g_{k,s} \leq \frac{g_{k+1}}{p^{2-s}}+1}$, so ${g_{k,s} \leq  \frac{q^{k}(q+1) - \sqrt{q}q^\frac{k+1}{2}}{p^{2-s}} +1}$ which gives the result since ${p^{2-s} \leq q^\frac{k}{2}}$ for all ${k\geq2}$.
\qed
\end{Proof}

\begin{Lemma}\label{lemme_delta}
 Let $q=p^2=4$. For all $k\geq 1$ and ${0 \leq s \leq 2}$, we set $D_{k,s}:=p^{s+1}q^{k-1}$. Then we have
\begin{enumerate}[i)]
	\item $\Delta g_{k,s} := g_{k,s+1} - g_{k,s} \geq D_{k,s}$,
	\item $N_1(H_{k,s}/\F_p) + 2N_2(H_{k,s}/\F_p) + 4N_4(H_{k,s}/\F_p) > 2D_{k,s}$.
\end{enumerate}
\end{Lemma}

\begin{Proof} \textit{i)} From Hurwitz Genus Formula, we know that \linebreak[4]${g_{k,s+1} - 1 \geq p(g_{k,s}-1)}$ for any integer ${k\geq1}$ and ${s=0,1}$, so \linebreak[4]${g_{k,s+1} - g_{k,s} \geq (p-1)(g_{k,s} -1)}$. Applying $s$ more times Hurwitz Genus Formula, we get ${g_{k,s+1} - g_{k,s} \geq (p-1)p^s(g_k -1)}$ thus for ${k\geq 4}$ we have ${g_{k,s+1} - g_{k,s} \geq (p-1)p^sq^k}$ because ${g_k > q^k}$ according to Lemma \ref{lemme_genre} i).

\textit{ii)} It is obvious since ${q^2-1 > p^2}$ and since  from Proposition~\ref{subfield} we have ${N_1(H_{k,s}/\F_2) + 2N_2(H_{k,s}/\F_2) + 4N_4(H_{k,s}/\F_2) \geq (q^2-1)q^{k-1}p^s}$.
\qed
\end{Proof}

\begin{Lemma}\label{lemme_bornesup}
Let $q=p^2=4$ and ${N_i(k,s) := N_i(H_{k,s}/\F_p)}$. For all ${k \geq 1}$ and ${s=0, 1, 2}$, we have
$$
\sup \Big \{ n \in \N \; | \;  2n \leq N_1(k,s) + 2N_2(k,s) + 4N_4(k,s) -2g_{k,s} -7 \Big \} \geq \frac{5}{2}q^{k-1}-\frac{7}{2}.
$$
\end{Lemma}

\begin{Proof}
From Proposition \ref{subfield} and Lemma \ref{lemme_genre} iii), we get 
\begin{eqnarray*}
N_1(k,s) + 2N_2(k,s) + 4N_4(k,s) -2g_{k,s} -7 & \geq & (q^2-1)q^{k-1}p^s\\
									& & \  - 2q^{k-1}(q+1)p^s-7 \\
						   & = & p^sq^{k-1}(q+1)(q-3) -7
\end{eqnarray*}
thus we have $\sup \Big \{ n \in \N \; | \;  2n \leq N_1(k,s) + 2N_2(k,s) + 4N_4(k,s) -2g_{k,s} -7 \Big \} \geq \frac{1}{2}p^sq^{k-1}(q+1)(q-3) - \frac{7}{2}$
and we get the result since ${q=4}$ and ${s\geq0}$.
\qed
\end{Proof}

\begin{Lemma}\label{lemme_placedegn}

Let $n$ be an integer ${\geq 2}$. Then there exists a step $H_{k,s}/\F_2$ of the tower $T_3/\F_2$ introduced in Section \ref{sectiondescent} such that both following conditions are verified:
\begin{enumerate}[(1)]
	\item there exists a place of degree $n$ in $H_{k,s}/\F_2$, 
	\item $N_1(H_{k,s}/\F_2)+2N_2(H_{k,s}/\F_2)+4N_4(H_{k,s}/\F_2) \geq 2n + 2g_{k,s}+7$.
\end{enumerate}
Moreover, the first step for which both conditions are verified is the first step for which (2) is verified.
\end{Lemma}

\begin{Proof} 
Let $q = p^2 = 4$. Fix $n\geq 28$. We first show that for all integers $k$ such that ${2 \leq k \leq \frac{1}{4}(n-12)}$, we have  ${2g_{k,s}+1 \leq p^\frac{n-1}{2}(p^\frac{1}{2}-1)}$ for any ${s \in \{0,1,2\}}$, so Condition (1) is verified according to Corollary 5.2.10 in \cite{stic}. Indeed for such an integer $k$, we have ${6 \leq \frac{n}{2} - 2k}$ i.e. ${p^6 \leq p^{\frac{n}{2}-2k}}$. Since ${5p^\frac{7}{2} \leq p^6}$, we get ${5p^\frac{7}{2} \leq p^{\frac{n}{2}-2k}}$ or equivalently ${5p^{2k+1} \leq p^{\frac{n-1}{2} -2}}$, which leads to ${5p^{2k+1} \leq p^\frac{n-1}{2}(p^\frac{1}{2}-1)}$. Now, let us show that ${2g_{k,s}+1 \leq 5p^{2k+1}}$. According to Lemma \ref{lemme_genre} iv), since ${k\geq2}$ we have for ${s=0,1,2}$:
\begin{eqnarray*}
	2g_{k,s}+1 & \leq & 2\frac{q^k(q+1)-q^\frac{k}{2}(q-1)}{p^{2-s}}+1\\
			  & = & 2\left(q^{k-1}(q+1)-q^\frac{k}{2}\frac{q-1}{q}\right)p^s+1\\
			  & = & 2q^{k-1}(q+1) p^s - 2q^\frac{k}{2}\frac{q-1}{q}p^s+1 \\
			  & \leq & 2q^{k-1}(q+1)p^s \ \ \mbox{ since } 2q^\frac{k}{2}\frac{q-1}{q}p^s \geq 1 \\
			  & = & 2p^{2(k-1)}(p^2+1)p^s \\
			  & = & 5p^{2k-1}p^s  \ \ \mbox{ since } p=2
\end{eqnarray*}
which gives the result since ${p^s \leq p^2}$. \\
We prove now that for ${k\geq \frac{1}{2}\log_p \left(\frac{4}{5}(2n+6)\right)}$, Condition (2) is verified. Indeed, for such an integer $k$, we have ${2n+6 \leq \frac{5}{4}p^{2k}}$, so ${2n+6 \leq \frac{5}{4}p^{2k}}p^s$ for ${s=0,1,2}$. Since ${p=2}$, we have ${\frac{5}{4}p^{2k} p^s= \big(p^4-1-p(p^2+1)\big)}p^{2k-2}p^s$, so we get 
\begin{equation}\label{ineqclef}
2n+p^{2k-1}(p^2+1)p^s + 6 \leq (p^4-1)p^{2k-2}p^s.
\end{equation}
 Recall that we got ${2g_{k,s} +1 \leq p^{2k-1}(p^2+1)}p^s$ in the first part of the proof, so ${2n+2g_{k,s}+7 \leq 2n +p^{2k-1}(p^2+1)p^s +6}$ and (\ref{ineqclef})  gives the result since we know from Proposition \ref{subfield} that \linebreak[4]${N_{1}(H_{k,s}/\F_p)+2N_{2}(H_{k,s}/\F_p) +4N_{4}(H_{k,s}/\F_p)\geq (q^2-1)q^{k-1}p^{s}}$.\\
 Finally, we have proved that for any integers ${n\geq 28}$ and ${k \geq 2}$ such that ${\frac{1}{2}\log_p \left(\frac{4}{5}(2n+6)\right) \leq k \leq \frac{1}{4}(n-12)}$, both Conditions (1) and (2) are verified. Note that for any ${n\geq 28}$, we have ${\frac{1}{2}\log_p \left(\frac{4}{5}(2n+6) \right)>2}$. Moreover the size of the interval ${\big[ \frac{1}{2}\log_p \left(\frac{4}{5}(2n+6)\right) ; \frac{1}{4}(n-12) \big]}$ is bigger than 1 as soon as ${n \geq 28}$, and this size increases with $n$. Hence, for any integer $n\geq28$, we know that there is an integer $k > 2$ in this interval and so there exists a corresponding step $H_{k,s}$. Moreover, the first step $H_{k,s}$, that is to say the smallest couple of integers ${(k,s)}$, for which both Conditions (1) and (2) are verified, is the first step for which Condition (2) is verified, since for all integers ${k \leq \frac{1}{4}(n-12)}$ there is a place of degree $n$ in  $H_{k,s}/\F_2$.
To conclude, we complete the proof by computing, for the firsts steps of the tower, the number of places of degree one, two, four and $n$ for ${n < 28}$. Using the KASH packages \cite{kash}, we obtain the following results :
\begin{enumerate}[a)]
	\item ${g(H_1/\F_2) = 0}$, ${N_1(H_1/\F_2) = 3}$, ${N_2(H_1/\F_2) = 1}$ and ${N_4(H_1/\F_2) = 3}$. Hence Condition (2) holds for all ${n \leq 5}$, and we check that moreover ${N_3(H_1/\F_2) >0}$ and ${N_5(H_1/\F_2) >0}$. So for any integer ${n \leq 5}$, the first step that verifies both Conditions (1) and (2) is $H_1/\F_2$. 
	\item  ${g(H_{1,1}/\F_2) = 2}$, ${N_1(H_{1,1}/\F_2) = 3}$, ${N_2(H_{1,1}/\F_2) = 1}$ and \linebreak[4]${N_4(H_{1,1}/\F_2) = 7}$. Hence Condition (2) holds for all ${n \leq 11}$, and we check that moreover ${N_i(H_{1,1}/\F_2) > 0}$ for all integers $i$ such that ${6 \leq i \leq 11}$. So for any integer $n$ such that ${6 \leq n \leq 11}$, the first step that verifies both Conditions (1) and (2) is $H_{1,1}/\F_2$. 
	\item ${g(H_2/\F_2) = 6}$, ${N_1(H_2/\F_2) = 3}$, ${N_2(H_2/\F_2) = 1}$ and ${N_4(H_2/\F_2) = 15}$. Hence Condition (2) holds for all ${n \leq 23}$, and we know that moreover ${N_i(H_2/\F_2) > 0}$ for all integers $i$ such that ${12 \leq i \leq 23}$ since we have ${2g(H_2/\F_2)+1 \leq 2^\frac{i-1}{2}(\sqrt{2}-1)}$. Indeed ${2g(H_2/\F_2)+1=13}$ and ${2^\frac{i-1}{2}(\sqrt{2}-1) \geq 2^\frac{12-1}{2}(\sqrt{2}-1) \geq 18}$ for all integers $i$ such that ${12 \leq i \leq 23}$. So for any integer $n$ such that ${12 \leq n \leq 23}$, the first step that verifies both Conditions (1) and (2) is $H_2/\F_2$. 
	\item ${g(H_{2,1}/\F_2) = 23}$, ${N_1(H_{2,1}/\F_2) = 4}$, ${N_2(H_{2,1}/\F_2) = 1}$ and \linebreak[4]${N_4(H_{2,1}/\F_2) = 28}$. Hence Condition (2) holds for all ${n \leq 32}$, and we know that moreover ${N_i(H_2/\F_2) > 0}$ for all integers $i$ such that ${24 \leq i \leq 27}$ since we have ${2g(H_{2,1}/\F_2)+1 \leq 2^\frac{n-1}{2}(\sqrt{2}-1)}$. Indeed ${2g(H_{2,1}/\F_2)+1=47}$ and ${2^\frac{i-1}{2}(\sqrt{2}-1) \geq 2^\frac{24-1}{2}(\sqrt{2}-1) \geq 1199}$ for all integers $i$ such that ${24 \leq i \leq 27}$. So for any integer $n$ such that ${24 \leq n \leq 27}$, the first step that verifies both Conditions (1) and (2) is $H_{2,1}/\F_2$. 
\end{enumerate}	
Note that, as in the first part of the proof, we have to use the step $(k,s+1)$ because Condition (2) is not verified for the step $(k,s)$.
\qed \\
\end{Proof}

Finally, we give the following lemma which ensures us that given a finite set of places $\mathcal P$ and a divisor $\D$, up to equivalence we can suppose that the support of $\D$ does not contain any place in $\mathcal P$.

\begin{Lemma}\label{lemme_depsupp}
Let $F/\F_q$ be an algebraic function field and \linebreak[4]${\mathcal P := \{P_1, \ldots, P_N\}}$ be a set of places of arbitrary degrees in $F/\F_q$. For any divisor $\D$, there exists a divisor $\D'$ such that $\D$ and $\D'$ are equivalents and ${ \mathcal P \cap \mathrm{supp} \, \D = \varnothing}$. 
\end{Lemma}

\begin{Proof}
Let us consider the integers ${n_1, \ldots, n_N}$ defined by ${n_i = 0}$ if ${P_i \notin \mathrm{supp} \, \D}$ and ${n_i = - \mathrm{ord}_{P_i} \, \D}$ if ${P_i \in \mathrm{supp} \, \D}$. According to Strong Approximation Theorem (cf \cite{stic}, Theorem 1.6.5), there exists an element ${x \in F/\F_q}$ such that for all integers ${i \in \{1, \ldots, N\}}$, ${v_{P_i}(x) = n_i}$ and for any place ${P \notin \mathcal P}$, ${v_P(x) \geq 0}$. Thus we have for all integers ${i \in \{1,\ldots, N\}}$, ${\mathrm{ord}_{P_i}\big(\D+(x)\big) = \mathrm{ord}_{P_i} \, \D + n_i = 0}$ i.e. the intersection ${\mathcal P \cap \mathrm{supp} \big( \D + (x) \big)}$ is empty, so ${\D' := \D + (x)}$ is a suitable \linebreak[4]$\D$-equivalent divisor.
\qed
\end{Proof}

\section{New bounds for the tensor rank}

\subsection{Adapted algorithm of type Chudnovsky and associated complexity}

In this section, we use places of degree one, two and four to obtain new results for the tensor rank of multiplication in any extension of the finite field $\F_2$.

First of all, we specialize the general algorithm presented in Theorem~\ref{theo_evalder} for places of degree one, two and four by using first derivative evaluations, i.e. with ${u_i \leq 2}$ for ${i=1, \ldots, N}$.

\begin{Proposition} \label{theo_deg4evalder}
Let \\
\vspace{.1em}
$\bullet$ $q$ be a prime power, \\
\vspace{.1em}
$\bullet$ $F/\F_q$ be an algebraic function field,\\
\vspace{.1em}
$\bullet$ $Q$ be a degree $n$ place of $F/\F_q$,\\ 
\vspace{.1em}
$\bullet$ $\D$ be a divisor of $F/\F_q$, \\
\vspace{.1em}
$\bullet$ ${\mathcal P}=\{P_1,\ldots,P_{N_1},P_{N_1+1}, \ldots, P_{N_1+N_2}, P_{N_1+N_2+1}, \ldots, P_{N_1+N_2+N_4}\}$ be \linebreak[4] \indent a set of $N_1$ places of degree one, $N_2$ places of degree two and \linebreak[4] \indent $N_4$ places of degree four.\\
\vspace{.1em}
$\bullet$ ${0 \leq l_1 \leq N_1}$, ${0 \leq l_2 \leq N_2}$ and ${0 \leq l_4 \leq N_4}$ be three integers.\\
We suppose that $Q$ and all the places in $\mathcal P$ are not in the support of $\D$ and that:
\begin{enumerate}[a)]
	\item the application
	$$
	Ev_Q: \Ld{} \rightarrow \F_{q^n}\simeq F_Q$$
	is onto,
	\item the application
	 $$
	 Ev_{\mathcal P}: \left \{
	\begin{array}{ccl}
   	\Ld{2} & \rightarrow & \F_q^{N_1} \times \F_q^{l_1}\times \F_{q^2}^{N_2} \times \F_{q^2}^{l_2} \times \F_{q^4}^{N_4} \times \F_{q^4}^{l_4}  \\
    	f & \mapsto & \big( f(P_1),\ldots,f(P_{N_1}), f'(P_1), \ldots, f'(P_{l_1}),f(P_{N_1+1}), \ldots,\\
			& & \ \  f(P_{N_1+N_2}), f'(P_{N_1+1}), \ldots, f'(P_{N_1+l_2}), f(P_{N_1+N_2+1}),\\
		& & \ \ \ldots, f(P_{N_1+N_2+N_4}), f'(P_{N_1+N_2+1}), \ldots, f'(P_{N_1+N_2+l_4}) \big)
	\end{array} \right .
	$$
is injective.
\end{enumerate}
Then  
$$
\mu_q(n)\leq  N_1+ 2l_1+ 3N_2+6l_2+\mu_q(4)\big(N_4 + 2l_4).
$$
\end{Proposition}

\begin{Proof}
Up to reindex the places, the result follows from Theorem \ref{theo_evalder} applied with $N=N_1+N_2+N_4$, ${\deg P_i=1}$ for ${i=1, \ldots, N_1}$, ${\deg P_i=2}$ for ${i=N_1+1, \ldots, N_1+N_2}$, ${\deg P_i=4}$ for ${i=N_1+N_2+1,}$ ${ \ldots, N}$ and 
$$u_i = \left \{
	\begin{array}{ll}
	2\mbox{,} & \mbox{if }1\leq i \leq l_1, \mbox{ or } N_1+1\leq i \leq N_1+l_2,\\
			& \ \  \mbox{ or } N_1+N_2 +1 \leq i \leq N_1+N_2+l_4 ,\\
	1\mbox{,} & \mbox{else.}
	\end{array} \right .$$
Recall that for all prime powers $q$, $\mu_q(2)=3$ and $\widehat{M_q}(2) \leq 3$.\\
Applying Theorem \ref{theo_evalder}, we get:
\begin{eqnarray*}
	\mu_q(n) & \leq & \sum_{i=1}^{l_1}\mu_q(1)\widehat{M_q}(2) + \sum_{i=l_1+1}^{N_1}\mu_q(1)\widehat{M_q}(1) + \sum_{i=N_1+1}^{N_1+l_2}\mu_q(2)\widehat{M_{q^2}}(2)\\
			&  &  \ + \sum_{i=N_1+l_2+1}^{N_1+N_2}\mu_q(2)\widehat{M_{q^2}}(1) + \sum_{i=N_1+N_2+1}^{N_1+N_2+l_4} \mu_q(4)\widehat{M_{q^4}}(2) \\
			 & & \ + \sum_{i=N_1+N_2+l_4+1}^{N} \mu_q(4) \widehat{M_{q^4}}(1)\\
			& \leq & 3l_1 + N_1 - l_1 +9l_2 + 3(N_2-l_2) + 3\mu_q(4) l_4 + \mu_q(4)(N_4-l_4)\\
			& = & N_1 +2l_1+3N_2 +6l_2 + \mu_q(4)(N_4+2l_4).
\end{eqnarray*}	
\qed
\end{Proof}\\

\Remark Note that if $l_1$, $l_2$ and $l_4$ are three integers such that the application $Ev_\mathcal{P}$ is injective, then for any other integers $L_1$, $L_2$ and $L_4$ such that ${l_1 \leq L_1 \leq N_1}$, ${l_2 \leq L_2 \leq N_2}$ and ${l_4 \leq L_4 \leq N_4}$ the injectivity of the application is still verified but we obtain a bigger bound for the bilinear complexity. Consequently, we will try to use the optimal integers $l_1$, $l_2$ and $l_4$, that is to say the smallest integers for which the application $Ev_\mathcal{P}$ is injective. In particular, if ${l_1=l_2=l_4=0}$ is a suitable choice, then we can multiply in $\F_{q^n}$ without using derivative evaluations. \\

\begin{Theorem}\label{result_deg4evalder}
Let $q$ be a prime power. Let $F/\F_q$ be an algebraic function field of genus $g$ and $N_i$ be a number of places of degree $i$ in $F/\F_q$. Let $l_1$, $l_2$, $l_4$ be three integers such that $0 \leq l_1 \leq N_1$, $0 \leq l_2 \leq N_2$ and $0 \leq l_4 \leq N_4$. If
\begin{enumerate}[i)]
	\item $N_n > 0$ (or $2g+1 \leq q^{\frac{n-1}{2}}(q^{\frac{1}{2}}-1)$),
	\item $N_1 +l_1 +2(N_2+l_2) +4(N_4+l_4) > 2n+2g+6$,
\end{enumerate}
then
$$
\mu_q(n) \leq \frac{\mu_q(4)}{2}\big(n+g+5\big) + \mu_q(4)l_4.
$$

In particular, 
$$
\mu_2(n) \leq \frac{9}{2}\big(n+g+5\big) +9l_4.
$$
\end{Theorem}

\begin{Proof} Let $Q$ be a place of degree $n$ in $F/\F_q$, which exists since~i). 
We can build a divisor $\D$ such that the application $Ev_Q$ defined previously is onto. Indeed from Corollary 3.4 in \cite{bariro}, there exists a zero dimension divisor $\D[R]$ of degree $g-5$. Let $\D$ be a divisor such that ${\D \thicksim \D[K] + Q - \D[R]}$, with $\D[K]$ a canonical divisor. According to Lemma \ref{lemme_depsupp}, we can choose $\D$ such that ${Q \notin \mathrm{supp} \, \D}$. Such a divisor $\D$ verifies  ${\deg \D = n+g+3}$ and by Riemann-Roch Theorem, we have  ${\dim(\D-Q) =4}$ since ${i(\D-Q)=\dim(\D[K]-\D+Q) = \dim \D[R] = 0}$. \linebreak[4]Moreover, by Riemann-Roch Theorem we get ${\dim \D \geq n+4}$. Consequently, $Ev_Q$ is onto since the dimension of its image verifies
$$
\dim \mathrm{Im}(Ev_Q) = \dim \D - \dim (\D-Q) \geq n.
$$
Let us set ${N:=N_1 +l_1 +2(N_2+l_2) +4(N_4+l_4)}$. According to ii), we know that ${N > 2n+2g+6}$ so w.l.o.g. we can assume that \linebreak[4]${N = 2n+2g+7+\epsilon}$ with ${\epsilon = 0,1,2,3}$. Let $\mathcal P$ be a set of $N_1$ places of degree one, $N_2$ places of degree two and $N_4$ places of degree four. According to Lemma \ref{lemme_depsupp}, we can suppose that no place in $\mathcal P$ is in the support of $\D$.Note that we can apply Proposition \ref{theo_deg4evalder} with the set of places  $\mathcal{P}$ by using $l_1$ derivative evaluations on places of degree one, $l_2$ derivative evaluations on places of degree two and $l_4$ derivative evaluations on places of degree four. Indeed, let us denote by $\D[A]$ the divisor \linebreak[4]$\D[A] := \sum_{i=1}^{N_1+N_2+N_4} P_i + \sum_{i=1}^{l_1} P_i + \sum_{i=1}^{l_2} P_{N_1+i} + \sum_{i=1}^{l_4} P_{N_1+N_2+i}$, then we have ${\deg \D[A] =N}$, so ${\deg (2 \D - \D[A]) <0}$ by ii) and ${\ker Ev_{\mathcal{P}} = \mathcal{L}\left(2\D - \D[A]\right)}$ is trivial. 
Thus, we get ${\mu_q(n) \leq N_1+2l_1+3N_2+ 6l_2+\mu_q(4)(N_4+2l_4)}$  by Proposition \ref{theo_deg4evalder}. Now  let us remark that this bound depends on the number of places of each degree we use in the second evaluation: the higher the degrees are, the bigger the bound is. Consequently, we must consider that ${N_1=N_2=0}$ corresponding to the \textit{worst case}. Then we obtain ${\mu_q(n) \leq \mu_q(4)(\frac{N}{4} + l_4)}$, which gives the result since \linebreak[4]${\frac{N}{4} \leq \frac{2n+2g+10}{4} = \frac{1}{2}(n+g+5)}$. In particular, for ${q=2}$ we get 
$$
\mu_2(n) \leq \frac{9}{2}\big(n+g+5\big) +9l_4
$$
since ${\mu_2(4)=9}$.

\qed
\end{Proof}\\

\subsection{Tensor rank in any extension of $\F_2$} Now we apply the results of the preceding section to the tower of Garcia-Stichtenoth $T_3/\F_2$ presented in Section \ref{sectiondescent}. We obtain two kinds of result: one which uses derivative evaluations and an other which does not. We will see later that we obtain a better bound for $M_2$ using derivative evaluations but this utilization is more complicated in practice and leads to an increase of linear complexity which can be inconvenient; so we present both technics. Moreover, although the best results are obtained using derivative evaluations, we still get an improvement of the best known bound for $M_2$ using simple evaluations.\\

\subsubsection{Bound for the tensor rank without using derivative evaluation} First of all, we apply the bound of Theorem \ref{result_deg4evalder} on the tower $T_3/\F_2$ with ${l_1=l_2=l_4=0}$.

\begin{Theorem}\label{bornesansevalder}
For any integer $n\geq2$, we have
$$
\mu_2(n) \leq \frac{45}{2}n + 85.5.
$$

\end{Theorem}

\begin{Proof}
Let $q=p^2=4$ and let us consider the sequence of algebraic function fields $T_3=\big\{H_{k,s}/\F_2\big\}$ introduced in Section \ref{sectiondescent}. We set ${M_{k,s} := N_1(H_{k,s}/\F_2)+2N_2(H_{k,s}/\F_2)+4N_4(H_{k,s}/\F_2)}$.
For any integer~$n$, we know by Lemma~\ref{lemme_placedegn} that there exists a step of the tower $T_3$ on which we can apply Theorem \ref{result_deg4evalder}. Let $H_{k,s}/\F_2$ be the first step of the tower that suits the hypothesis of Theorem \ref{result_deg4evalder} with ${l_1=l_2=l_4=0}$. According to Lemma \ref{lemme_placedegn}, this step is determined by the smallest integers $k$ and $s$ such that ${2n \leq M_{k,s}-2g_{k,s}-7}$, so ${2n > M_{k,s-1}-2g_{k,s-1}-7}$.
For any integer ${k \geq 1}$ and for any integer ${s=0,1,2}$, we have ${g_{k,s}\leq q^{k-1}(q+1)p^{s}}$ by Lemma \ref{lemme_genre} iii). Moreover, since ${M_{k,s-1} \geq (q^2-1)q^{k-1}p^{s-1}}$ by Proposition \ref{subfield}, we \linebreak[4]obtain ${2n > (q^2-2q-3)q^{k-1} p^{s-1}-7}$. Then since ${q=4}$, we have\linebreak[4] ${q^2-2q+3 = (q+1)(q-3)=q+1}$, which leads to \linebreak[4] ${2np > (q+1)q^{k-1}p^s - 7p \geq g_{k,s}-7p}$ and it follows that \linebreak[4]${g_{k,s} \leq 2np+7p}$, so

$$
{\mu_{2}(n) \leq \frac{9}{2}\big(n+g_{k,s}+5\big) \leq \frac{9}{2}n\left(1+2p\right) + \frac{9}{2}(7p+5)}
$$
by Theorem \ref{result_deg4evalder}, which gives the result since $p=2$.
\qed
\end{Proof}

\subsubsection{Bound for the tensor rank using derivative evaluations} Here, we apply results of Theorem \ref{result_deg4evalder} with an optimal number of derivative evaluations.

\begin{Theorem}\label{borneevalder}
For any integer $n\geq2$, we have
$$
\mu_2(n) \leq \frac{477}{26}n +\frac{45}{2}.
$$

\end{Theorem}

\begin{Proof}
For any integer $n$, we know by Lemma~\ref{lemme_placedegn} that there exists a step of the tower $T_3/\F_2$ on which we can apply Theorem~\ref{result_deg4evalder} with ${l_1=l_2=l_4=0}$. We set ${M_{k,s}:= N_1(H_{k,s}/\F_2)+2N_2(H_{k,s}/\F_2) +}$ ${4N_4(H_{k,s}/\F_2)}$ for any step $H_{k,s}/\F_2$, with ${k \geq 0}$ and ${s=0,1}$. Let \linebreak[4]$H_{k,s+1}/\F_2$ be the first step of the tower that suits the hypothesis of Theorem \ref{result_deg4evalder} with ${l_1=l_2=l_4=0}$ i.e. $k$ and $s$ are integers such that ${M_{k,s+1} > 2n + 2g_{k,s+1} +6}$ and ${M_{k,s} \leq 2n + 2g_{k,s} +6}$. We denote by $n_0^{k,s}$ the biggest integer such that ${M_{k,s} > 2n_0^{k,s} + 2g_{k,s} +6}$ i.e. \linebreak[4]${n_0^{k,s} := \sup \Big \{ n \in \N \; | \;  2n \leq M_{k,s} -2g_{k,s} -7 \Big \}}$. To multiply in $\F_{2^n}$, we have the following alternative:
\begin{enumerate}[a)]
	\item to use the algorithm on the step $H_{k,s+1}$. In this case, a bound for the bilinear complexity is given by Theorem \ref{result_deg4evalder} applied with ${l_1=l_2=l_4=0}$:
	$$
	\mu_2(n) \leq \frac{9}{2}(n+g_{k,s+1}+5) = \frac{9}{2}(n_0^{k,s}+g_{k,s}+5) +\frac{9}{2}(n-n_0^{k,s}+ \Delta g_{k,s}).
	$$
	Recall that ${\Delta g_{k,s}:= g_{k,s+1}-g_{k,s}}$.
	\item to use the algorithm on the step $H_{k,s}$ with derivative evaluations on $l_1$ places of degree one, $l_2$ places of degree two and $l_4$ places of degree four, where $l_i$ satisfies ${l_i \leq N_i(H_{k,s}/\F_2)}$ for ${i=1,2,4}$ and \linebreak[4]${M_{k,s}+l_1+2l_2+4l_4 > 2n + 2g_{k,s} + 6}$. One can check that this condition is verified as soon as ${l_1+2l_2+4l_4\geq 2(n-n_0^{k,s})}$,  so Theorem~\ref{result_deg4evalder} gives ${\mu_2(n) \leq \frac{9}{2}\big(n+g\big) +9l_4 + \frac{45}{2}}$. W.l.o.g, we can suppose that ${l_1+2l_2+4l_4= 2(n-n_0^{k,s})+ \epsilon}$ with ${\epsilon =0,1,2,3}$. Moreover, we must consider that ${l_1=l_2=0}$ corresponding to the \textit{worst case}. Thus we have ${4l_4= 2(n-n_0^{k,s})+ \epsilon \leq  \left[\frac{1}{2}(n-n_0^{k,s})\right]+1}$ with ${\left[ \cdot \right]}$ denoting the floor function, and we obtain the following bound for the bilinear complexity: 
\begin{eqnarray*}
	\mu_2(n) & \leq & \frac{9}{2}(n+g_{k,s}+5) + 9\left(\left[\frac{1}{2}(n-n_0^{k,s})\right] + 1\right) \\
			& \leq & \frac{9}{2}(n_0^{k,s}+g_{k,s}+5) +9(n-n_0^{k,s}+1).
\end{eqnarray*}
\end{enumerate}
Thus, if the integers $l_i$ such that ${l_1+2l_2+4l_4= 2(n-n_0^{k,s})+ \epsilon}$ with ${\epsilon =0,1,2,3}$, verifies ${l_i \leq N_i(H_{k,s}/\F_2)}$ for ${i=1,2,4}$, i.e. \linebreak[4]${2(n-n_0^{k,s})+ \epsilon \leq M_{k,s}}$ then case b) gives a better bound as soon as ${\Delta g_{k,s} > n-n_0^{k,s} +2}$.

For $x\in \R^+$ such that ${M_{k,s+1} > 2\left[x\right] + 2g_{k,s+1} +6}$ and \linebreak[4]${M_{k,s} \leq 2\left[ x \right] + 2g_{k,s} +6}$, we define the function $\Phi_{k,s}(x)$ as follow:
$$
\Phi_{k,s}(x) = \left\{\begin{array}{ll}
				9(x-n_0^{k,s}) + \frac{9}{2}(n_0^{k,s}+g_{k,s}+5) +9 & \mbox{if } x- n_0^{k,s} +2 < D_{k,s}\\
				\frac{9}{2}(x-n_0^{k,s}) + \frac{9}{2}(n_0^{k,s}+g_{k,s}+5+\Delta g_{k,s}) & \mbox{else}.
			      \end{array} \right.
$$
Recall that ${D_{k,s}}$ was defined in Lemma \ref{lemme_delta} as ${p^{s+1}q^{k-1}}$.\\
Note that if  ${x- n_0^{k,s} +2 < D_{k,s}}$, then according to Lemma \ref{lemme_delta} we have both 
$${x- n_0^{k,s} +2 < \Delta g_{k,s},}$$
 so case b) gives a better bound for the bilinear complexity,
and
$${2(x- n_0^{k,s}) + \epsilon < 2D_{k,s} \leq M_{k,s}} \mbox{ for } \epsilon =0,1,2,3,$$
 so we can proceed as in case b) since there are enough places of each degree to use derivative evaluations on $l_1$ places of degree one, $l_2$ \linebreak[4] places of degree two and $l_4$ places of degree four with \linebreak[4]${l_1+2l_2+4l_4= 2(n-n_0^{k,s})+ \epsilon}$.

We define the function $\Phi$ for all ${x\geq0}$ as the minimum of the functions $\Phi_{k,s}$ for which $x$ is in the domain of $\Phi_{k,s}$. This function is piecewise linear with two kinds of piece: those which have slope $\frac{9}{2}$ and those which have slope $9$. Moreover, since the y-intercept of each piece grows with $k$ and $s$, the graph of the function $\Phi$ lies below any straight line that lies above all the points ${\big(n_0^{k,s}+D_{k,s} -2, \Phi(n_0^{k,s}+D_{k,s} -2)\big)}$, since its are the \textit{vertices} of the graph. Let ${X:= x_0^{k,s} +D_{k,s} - 2}$, then
$$
\Phi(X) = \frac{9}{2}(X + g_{k,s+1} +5) = \frac{9}{2}\left(1+ \frac{g_{k,s+1}}{X}\right)X +\frac{45}{2}.
$$
We want to give a bound for $\Phi(X)$ which is independent of $k$ and $s$.

Lemmas \ref{lemme_genre} iii) and \ref{lemme_bornesup} give
\begin{eqnarray*}
	\frac{g_{k,s+1}}{X} & \leq & \frac{q^{k-1}(q+1)p^{s+1}}{\frac{5}{2}q^{k-1}-\frac{7}{2} + p^{s+1}q^{k-1} - 2}\\
				      & = & \frac{q+1}{\frac{5}{2p^{s+1}} +1- \frac{11}{2q^{k-1}p^{s+1}}}\\
				      & \leq & \frac{5}{\frac{13}{8} - \frac{11}{4\cdot4^{k-1}}}\\
				      & \leq & \frac{40}{13}.
\end{eqnarray*}
Thus, the graph of the function $\Phi$ lies below the line ${y=\frac{9}{2}\left(1+\frac{40}{13}\right)x+\frac{45}{2}}$. In particular, we get
$$
\Phi(n) \leq \frac{9}{2}\left(1+\frac{40}{13}\right)n+\frac{45}{2}.
$$ \qed
\end{Proof}

\section{New asymptotical bounds for the tensor rank}

Without using derivative evaluation, we obtain from Theorem \ref{bornesansevalder} the following bound for~$M_2$:

$$
M_2 \leq 22.5,
$$
 which is better than the best known bound recalled in Proposition \ref{theo_M2_2285}.

However, by using derivative evaluations, we obtain a better bound for $M_2$, namely:

\begin{Theorem}\label{theo_M2_2086}
$$
M_2 \leq \frac{477}{26} \approx 18.35.
$$
\end{Theorem}

\begin{Proof} 
It follows from Theorem \ref{borneevalder}.
\qed
\end{Proof}

%

\end{document}